\documentclass[a4paper,12pt]{article}

\newtheorem{proposition}{PROPOSITION}[section]
\newtheorem{theorem}[proposition]{THEOREM}

\usepackage{latexsym, amsfonts}

\newcommand{\PP}{\mbox{$\mathbb P$}}

\newcommand{\call}[1]{{\mathcal{ #1}}}

\newcommand{\rist}[1]{\mbox{\raisebox{-1ex}{$|#1$}}}
\newcommand{\ml}{linear morphism}
\newcommand{\mli}{linear morphisms}
\newcommand{\applicaz}[3]{\mbox{$#1\,:\,#2
  \rightarrow#3$}}

\newcommand{\lsb}{\mbox{\boldmath$[$}}
\newcommand{\rsb}{\mbox{\boldmath$]$}}
\newcommand{\ovf}{\overline f}
\newcommand{\ovfx}[1]{\overline{{}_{#1}\!f}}

\begin{document}

\title{On embedded products of Grassmannians%
\thanks{This work was performed
in the context of the fourth protocol of
  scientific and technological cooperation between Italy and Austria
  (plan n.~10).
  Partial support was provided by the project ``Strutture Geometriche,
  Combinatoria e loro Applicazioni'' of M.U.R.S.T.}}
\author{Hans Havlicek\\Technische Universit\"at Wien,\\
Abteilung f\"ur Lineare Algebra und Geometrie,\\ Wiedner Hauptstrasse
8-10/1133,\\
A-1040 Wien, Austria\\havlicek@geometrie.tuwien.ac.at \and Corrado Zanella\\
Dipartimento di Matematica Pura ed Applicata,\\ via Belzoni 7,\\ I-35131
Padova,\\ Italy\\zanella@math.unipd.it}

\date{}

\maketitle

\begin{abstract}
  Let $\Gamma'$ and $\Gamma$ be two Grassmannians.
  The standard embedding \applicaz{\phi}{\Gamma'\times\Gamma}{\overline{\PP}}\
  is obtained by combining the Pl\"ucker and Segre embeddings.
  Given a further embedding
  \applicaz{\eta}{\Gamma'\times\Gamma}{\PP'},
  we find a sufficient condition for the existence of
  $\alpha\in{\rm Aut}(\Gamma)$ and of a collineation
  \applicaz{\psi}{\overline{\PP}}{\PP'}\ such that
  $\eta=({\rm id}_{\Gamma'}\times\alpha)\phi\psi$.

\noindent Keywords: Segre variety; product space; Grassmann variety; projective
embedding

\noindent A.M.S.\ classification number: 51M35.
\end{abstract}

\sloppy


\section{Introduction}

\subsection{Background}
Several authors have proved that the
classical embeddings of geometries such as Grassmannians and
product spaces are essentially unique. For example,
Havlicek~\cite{Hav81} showed that every embedding of a Grassmann
space can be represented as the product of the standard embedding,
which is obtained by means of Pl\"ucker coordinates, and a \ml\
between projective spaces (essentially a projection between
complementary subspaces).

Such a strong universal property does not hold in general for
product spaces. However, if
\applicaz{\gamma}{\PP_1\times\PP_2}{\overline{\PP}}\ is the Segre
embedding, and \applicaz{\chi}{\PP_1\times\PP_2}{\PP'}\ is any
embedding, then there exist $\alpha\in{\rm Aut}(\PP_2)$ and a \ml\
\applicaz{\psi}{\overline{\PP}}{\PP'}\ such that $\chi=({\rm
id}_{\displaystyle\PP_1}\times\alpha)\gamma\psi$ (cf.\
Zanella~\cite{Zan96}).

As a consequence of the results in \cite{Hav81,Zan96}, the image
of any embedding of a Grassmann space or a product space is
projectively equivalent to a projection of the related variety.
So, the incidence geometrical characterizations given by Tallini
and several other authors are also intrinsic characterizations of
those varieties (cf.\ the surveys \cite{BiZa99,Tal86}).

We are attempting to give a result analogous to
\cite{Zan96} for an embedding $\eta$ of the product of two
Grassmannians $\Gamma'$ and $\Gamma$. As a first step, in this paper we
characterize the product of two Grassmannians up to collineations.

\subsection{Preliminaries}
A {\em semilinear space\/} is a
pair $\Sigma=(\call P,\call G)$, where $\call P$ is a set, whose
elements are called {\em points\/}, and $\call G\subseteq2^{\call
P}$. The elements of $\call G$ are {\em lines\/}. The axioms
defining a semilinear space are the following: ({\em
i\/})~$|g|\geq2$ for every line $g$; ({\em
ii\/})~$\displaystyle\bigcup_{g\in\call G}g=\call P$; ({\em
iii\/})~$g,h\in\call G$, $g\neq h$ $\Rightarrow$ $|g\cap h|\leq1$.
Two points $X,Y\in\call P$ are {\em collinear\/}, $X\sim Y$, if a
line $g$ exists such that $X,Y\in g$ (for $X\neq Y$ we will also
write $XY:=g$). An {\em isomorphism\/} between the semilinear
spaces $(\call P,\call G)$ and $(\call P',\call G')$ is a
bijection \applicaz{\alpha}{\call P}{\call P'}\ such that both
$\alpha$ and $\alpha^{-1}$ map lines onto lines.

The {\em join\/} of $\call M_1,\call M_2\subset \call P$ is:
\[
  \call M_1\vee\call M_2:=\call M_1\cup\call M_2\cup
  \bigg(
    \bigcup_{
        X_i\in\call M_i\atop X_1\sim X_2,\,X_1\neq X_2
    }X_1X_2
  \bigg).
\]

If $X$ is a point, we will often write $X$ instead of $\{X\}$.

Let $S$ and $T$ be sets. A {\em generalized mapping\/}, briefly
{\em g-map\/}, \applicaz fST\ is a mapping of a subset ${\bf
D}(f)$ of $S$ into $T$. ${\bf D}(f)$ is the {\em domain\/} of $f$
and ${\bf A}(f):=S\setminus{\bf D}(f)$ is the {\em exceptional
set\/}. If $X\in{\bf A}(f)$, then $Xf=\emptyset$. If ${\bf
D}(f)=S$, then $f$ is called a {\em global\/} g-map.

Let $\PP'=(\call P',\call G')$ be a projective space. A {\em
\ml\/} \applicaz{\chi}{\Sigma}{\PP'}\ is a g-map of $\call P$ into
$\call P'$ satisfying the following axioms (L1) and
(L2)~\cite{Bra73,Hav81}:

(L1)\hspace{1cm}$(X\vee Y)\chi=X\chi\vee Y\chi$ for $X,Y\in\call
P$, $X\sim Y$;

(L2)\hspace{1cm}$X,Y\in{\bf D}(\chi)$, $X\chi=Y\chi$,  $X\neq Y$,
$X\sim Y$ $\Rightarrow$ $\exists A\in XY$ such that
$A\chi=\emptyset$.

The \ml\ $\chi$ is called {\em embedding\/} if it is global and
injective. It should be noted the last definition is somewhat
particular, since for instance the inclusion of an affine space
into its projective extension is not an embedding.

The {\em (projective) rank} of $\chi$, rk$\chi$, is the projective
dimension of $\lsb \call P\chi\rsb$, where the square brackets
\lsb$\,$\rsb\ denote projective closure.


\section{Bilinear g-maps}

Let $\Sigma'$, $\Sigma''$ be semilinear spaces, $\PP$ a projective
space, and \applicaz f{\Sigma'\times\Sigma''}{\PP} a g-map. If for
every point $P$ of $\Sigma'$ the g-map
\applicaz{{}_P\!f}{\Sigma''}{\PP}\ defined by $X{}_P\!f:=(P,X)f$
is a \ml, then we say that $f$ is {\em right linear\/}. The
definition of a left linear g-map is similar. If the g-map $f$ is
both left linear and right linear, then it is called {\em
bilinear\/}. The bilinear mappings are exactly the \mli\ of the
product spaces. If $f$ is a bilinear g-map and for every point $P$
of $\Sigma'$ the g-map ${}_P\!f$ is an embedding, then we say that
$f$ is a {\em right embedding\/}. So, a right embedding is a
special type of global \ml.

Let $\Sigma$ be a semilinear space embedded in an $n$-dimensional
projective space \PP\ (that is, the inclusion is an embedding).
This semilinear space satisfies the {\em chain condition\/} (with
respect to \PP) if there are a plane $\call E$ and $n-2$ lines
(say $\ell_3$, $\ell_4$, $\ldots$, $\ell_n$) of $\Sigma$, such
that for every $i=3,\ldots,n$, $\dim\lsb\call
E\cup\ell_3\cup\ell_4\cup\ldots\cup\ell_i\rsb=i$.

Let $\Sigma$ be a semilinear space embedded in a projective space
\PP. If every embedding \applicaz f{\Sigma}{\PP'}\ can be uniquely
extended to a \ml\ \applicaz {\overline f}{\PP}{\PP'}, then we say
that $\Sigma$ is {\em universally embedded\/} in \PP. Obviously,
if $\Sigma$ is universally embedded in \PP, then $\Sigma$ spans
\PP. For instance: Every Grassmann variety is universally embedded
in its ambient space~\cite{Hav81}; a set of three pairwise skew
lines of a projective space \PP\ of dimension 3 is not universally
embedded in \PP.

\begin{proposition}\label{p4}
  Let $\Sigma'$ and $\Sigma''$ be two semilinear spaces, and\\
  \applicaz F{\Sigma'\times\Sigma''}{\PP'}\ a right embedding.
  Assume that $\Sigma''$ is universally embedded in a projective space $\PP$
  of dimension $n$.
  Let $\ell$ be a line of $\Sigma'$ and $\ell_1$, $\ell_2$ lines of $\Sigma''$
  such that $\ell_1\cap\ell_2$ is a point $P^*$.
  If
  \begin{equation}\label{ep4}
    \dim\lsb(\ell\times\Sigma'')F\rsb\geq 2n+1,
  \end{equation}
  then for $i=1,2$, $(\ell\times\ell_i)F$ is a hyperbolic quadric of the
  threedimensional subspace $U_i=\lsb(\ell\times\ell_i)F\rsb$ of $\PP'$.
  Furthermore, $U_1\cap U_2=(\ell\times P^*)F$.
\end{proposition}

{\em Proof\/}. Let $i\in\{1,2\}$. In order to prove that
$(\ell\times\ell_i)F$ is a hyperbolic quadric, it is enough to
check that given two distinct points $A,B\in\ell$, the lines
$(A\times\ell_i)F$ and $(B\times\ell_i)F$ are skew. Since the
embedding \applicaz{{}_A\!F}{\Sigma''}{\PP'}, mapping $X$ into
$(A,X)F$, can be linearly extended to $\PP$, we have
$\dim\lsb(A\times\Sigma'')F\rsb\leq n$. Then, in view of (\ref{ep4}),
we obtain:
\begin{equation}\label{referee}
  \mbox{For every }A,B\in\ell,\ A\neq B:\
  \lsb(A\times\Sigma'')F\rsb\cap\lsb(B\times\Sigma'')F\rsb=\emptyset.
\end{equation}
So, $(\ell\times\ell_i)F$ is a hyperbolic quadric.
As a further consequence of (\ref{referee}), since
$(A,P^*)F\neq(B,P^*)F$, and $F$ is left linear, $(\ell\times
P^*)F$ is a line.
Such a line
is contained in $U_1\cap U_2$.
If $\dim(U_1\cap U_2)>1$, then $U_1\cap U_2$ contains
a plane $\tilde{\call E}$ that is tangent to both quadrics
$(\ell\times\ell_i)F$, $i=1,2$. Then $\tilde{\call E}$ contains
two lines of type $(Q_i\times\ell_i)F$, with $Q_i\in\ell$,
$i=1,2$.

If $Q_1=Q_2$, then $(Q_1\times\ell_1)F\neq(Q_2\times\ell_2)F$
since we have a right embedding, whence
$\tilde{\call
E}\subset\lsb(Q_1\times\Sigma'')F\rsb$; but $\tilde{\call E}$ also
contains a point of type
$(Q^*,P^*)F\in\lsb(Q^*\times\Sigma'')F\rsb$ with
$Q^*\in\ell\setminus Q_1$. This implies
$\lsb(Q_1\times\Sigma'')F\rsb\cap\lsb
Q^*\times\Sigma'')F\rsb\neq\emptyset$, contradicting (\ref{referee}).

If $Q_1\neq Q_2$ we can obtain a similar contradiction because
$(Q_1\times\ell_1)F$ and $(Q_2\times\ell_2)F$ have a common point,
and
$(Q_1\times\Sigma'')F\cap(Q_2\times\Sigma'')F\neq\emptyset$.$\Box$

Let $\Sigma'$ and $\Sigma$ be two semilinear spaces, and \applicaz
f{\Sigma'\times\Sigma}{\PP'}\ a right embedding. If $\Sigma$ is
universally embedded in a projective space \PP, then for every
point $P$ of $\Sigma'$ the g-map \applicaz{{}_P\!f}{\Sigma}{\PP'}\
has a unique linear extension
\applicaz{\overline{{}_P\!f}}{\PP}{\PP'}. So, by setting
$(P,Q)\overline f:=Q\overline{{}_P\!f}$, we obtain a right linear
g-map \applicaz{\overline f}{\Sigma'\times\PP}{\PP'}\ which
extends $f$.

\begin{proposition}\label{p5}
  Let $\Sigma'$ and $\Sigma$ be two semilinear spaces.
  Assume that
  (i)~$\Sigma$ is universally embedded in an $n$-dimensional projective
  space \PP\ and satisfies the chain condition;
  (ii)~\applicaz f{\Sigma'\times\Sigma}{\PP'}\ is a right embedding;
  (iii)~for every line $\ell$ of $\Sigma'$,
  $\dim\lsb(\ell\times\Sigma)f\rsb\geq 2n+1$.

  Then
  the right linear extension
  \applicaz{\overline f}{\Sigma'\times\PP}{\PP'}\ is bilinear.
\end{proposition}

{\em Proof\/}. {\em (a)}\hspace{1em} Taking into account the chain
condition for $\Sigma$ we define
\[
  S_1:=\emptyset, \hspace{1em} S_2:=\call E, \hspace{1em}
  S_i:=\call E\cup\ell_3\cup\ell_4\cup\ldots\cup\ell_i,\ i=3,4,\ldots,n,
\]
\[
  T_i:=\lsb S_i\rsb,\ i=1,2,\ldots,n.
\]
Set $(\call P,\call G):=\Sigma$ and let $\call G_i$ be the line
set of $T_i$ ($i=1,2,\ldots,n$). We will show that the semilinear
spaces
\[
  \Sigma_i:=(\call P\cup T_i,\call G\cup\call G_i)
\]
are universally embedded in \PP\ for all $i=1,2,\ldots,n$. So let
\applicaz{F_i}{\Sigma_i}{\tilde{\PP}}\  be an embedding in some
projective space $\tilde{\PP}$. By ({\em i\/}), the embedding
$F_i\rist{\Sigma}$ can be extended uniquely to a \ml\ \applicaz
{F'_i}{\PP}{\tilde{\PP}}. Furthermore, $F'_i\rist{T_i}$ and
$F_i\rist{T_i}$ are two \mli\ which agree on $S_i$. In particular,
$F_i\rist{\call E}=F'_i\rist{\call E}$ is a collineation.
From \cite{Hav81}, Satz 1.3, $F_i\rist{\ell_j}=F'_i\rist{\ell_j}$ for
$j=3,4,\ldots,i$ and an easy induction on $j$, we obtain that
$F'_i\rist{T_i}=F_i\rist{T_i}$; so, $\Sigma_i$ is universally
embedded in \PP, as required.

{\em (b)}\hspace{1em} For a point $P$ of \PP, let \applicaz
{\overline {f}_P}{\Sigma'}{\PP'}\ be the g-map defined by
$X\overline {f}_P:=(X,P)\overline f$. It is enough to prove the
following statement by induction on $i=2,3,\ldots,n$:

\vspace{2mm}

(P$_i$) {\em If $P\in T_i\setminus (T_{i-1}\cup\call P)$, then
$\overline {f}_P$ is a \ml.\/}

\vspace{2mm}

First, (P$_2$) is trivial. Next, assume that (P$_i$) holds for
some $1<i<n$. Take a point $P\in T_{i+1}\setminus (T_i\cup\call
P)$, and define $g:=(\ell_{i+1}\vee P)\cap T_i$.

Let $\ell$ be any line of $\Sigma'$; we shall prove that
$\overline {f}_P\rist{\ell}$ is injective and that $\ell\overline
{f}_P$ is a line of $\PP'$; this will conclude the proof. For any
point $A$ of $\Sigma'$ take into account the \ml\
\applicaz{{\overline{{}_A\!f}}}{\PP}{\PP'}. There is a line $a$ of
$\Sigma'$ with $A\in a$ and a point $A'\in a\setminus A$. Then the
subspace spanned by $(a\times\Sigma)f$ is also spanned by the
image of ${\overline{{}_A\!f}}$ together with the image of
${\overline{{}_{A'}\!f}}$, so that ({\em iii\/}) implies
rk$({\overline{{}_A\!f}})=n$. Hence ${\overline{{}_A\!f}}$ is an
embedding and
\[
  F:=\overline f\rist{\Sigma'\times(T_i\cup\call P)}
\]
is a right embedding which allows to apply prop.~\ref{p4} with
$\Sigma'':=\Sigma_i$, $\ell_1:=g$, $\ell_2:=\ell_{i+1}$,
$P^*:=g\cap\ell_{i+1}$. Let $r'$ and $r''$ be two lines through
$P$ such that $r'\vee r''=\ell_1\vee\ell_2$,
$\ell_1\cap\ell_2\not\in r'\cup r''$. Furthermore, let
$B'_j:=r'\cap\ell_j$ and $B''_j:=r''\cap\ell_j$ ($j=1,2$).
From prop.~\ref{p4} the five lines
\begin{equation}\label{hans}
  (\ell\times P^*)f,\ (\ell\times B'_1)\ovf,\ (\ell\times B'_2)f,\
  (\ell\times B''_1)\ovf,\ (\ell\times B''_2)f
\end{equation}
are mutually skew and their span is 5-dimensional. As $L$ varies
in $\ell$, the four mappings $(L,P^*)f\mapsto(L,B'_j)\ovf$,
$(L,P^*)f\mapsto(L,B''_j)\ovf$ ($j=1,2$) are projectivities. Hence
$(L,B'_1)\ovf\mapsto(L,B'_2)f$ is a projectivity too and the
family of lines $(L\times r')\ovf$ with $L\in\ell$ is a regulus
with transversal lines $(\ell\times B'_1)\ovf$ and $(\ell\times
B'_2)f$. Fix one $L\in\ell$: By the linearity of $\ovfx{L}$, the
line $(L\times r')\ovf$ carries the point $(L,P)\ovf$. Since the
lines of a regulus are mutually skew, the mapping $(\ell\times
P)\ovf\ \rightarrow\ \{(L\times r')\ovf|L\in\ell\}$,\hspace{1em}
$(L,P)\ovf\ \mapsto\ (L,r')\ovf$ is bijective. Hence
$\ovf_P\rist{\ell}$ is injective. Similar arguments hold true for
$r''$. Now (\ref{hans}) implies that
\[
  ((\ell\times B'_1)\ovf\vee(\ell\times B'_2)f)\cap
  ((\ell\times B''_1)\ovf\vee(\ell\times B''_2)f)
\]
is a line $\tilde{\ell}$ which contains $\ell\ovf_{P}$. Hence
$\tilde{\ell}$ is a common transversal line of the two reguli
$\{(L\times r'){\ovf}|L\in\ell\}$ and $\{(L\times
r''){\ovf}|L\in\ell\}$ so that $\tilde{\ell}=\ell\ovf_{P}$.$\Box$

Let $F$ be a commutative field. We consider the integers
$0<h'<N'$, $0<h<N$, $n'={N'+1\choose h'+1}-1$,  $n={N+1\choose
h+1}-1$, the Grassmann spaces $\Gamma'=\Gamma^{h'}(\PP_{N',F})$,
$\Gamma=\Gamma^{h}(\PP_{N,F})$, with their Pl\"ucker embeddings
\applicaz{\wp'}{\Gamma'}{\PP_{n',F}},
\applicaz{\wp}{\Gamma}{\PP_{n,F}}. Let
\applicaz{\gamma}{\PP_{n',F}\times\PP_{n,F}}{\overline{\PP}}\ be
the Segre embedding ($\dim\overline{\PP}=(n'+1)(n+1)-1$). The {\em
standard embedding\/} of $\Gamma'\times\Gamma$ is the composition
$\phi:=(\wp'\times\wp)\gamma$.

\begin{theorem}\label{p6}
  Let \applicaz{\eta}{\Gamma'\times\Gamma}{\PP'}\ be an embedding such that
  ${\rm rk}\eta\geq{\rm rk}\phi$ and $\lsb{\rm im}\eta\rsb=\PP'$.
  Then there are $\alpha\in{\rm Aut}(\Gamma)$ and a collineation
  \applicaz{\psi}{\overline{\PP}}{\PP'}\ such that
  $\eta=({\rm id}_{\Gamma'}\times\alpha)\phi\psi$.
  Furthermore, ${\rm im}\eta$ is projectively equivalent to the image
  of the standard embedding
  of $\Gamma'\times\Gamma$.
\end{theorem}

{\em Proof\/}. We will identify $\Gamma'$ and $\Gamma$ with their
images under $\wp'$ and $\wp$, respectively. Let $A$ be a point of
$\Gamma'$, and \applicaz{{}_A\!\eta}{\Gamma}{\PP'}\ the embedding
of $\Gamma$ defined by $X{}_A\!\eta:=(A,X)\eta$. By the main
result in \cite{Hav81}, such an embedding can be extended to a
\ml\ \applicaz{{}_A\!\eta'}{\PP_{n,F}}{\PP'}. So, by setting
$(A,B)\overline{\eta}:=B{}_A\!\eta'$, we have a right linear g-map
\applicaz {\overline{\eta}}{\Gamma'\times\PP_{n,F}}{\PP'}. By
prop.~\ref{p5}, $\overline{\eta}$ is bilinear. By \cite{Hav81}
again, $\overline{\eta}$ has a left linear extension
\applicaz{{\overline{\overline{\eta}}}}{\PP_{n',F}\times\PP_{n,F}}{\PP'}.
We can apply the symmetric of prop.~\ref{p5}; condition
$\dim\lsb(\Gamma'\times\ell)\rsb\geq 2n'+1$ is a consequence of
the assumption ${\rm rk}\eta\geq{\rm rk}\phi$. By the main theorem
in \cite{Zan96}, there are $\tilde{\alpha}\in{\rm Aut}(\PP_{n,F})$
and a collineation \applicaz{\tilde\psi}{\overline{\PP}}{\PP'},
such that ${\overline{\overline{\eta}}}=({\rm
id}_{\displaystyle\PP_{n',F}}\times\tilde{\alpha})\gamma\tilde{\psi}$.
Every collineation of a projective space transforms both Grassmann
and Segre varieties in projectively equivalent ones (cf.\
e.g.~\cite{Zan95} (1.1)). Then
$\tilde\alpha=\tilde\alpha_1\tilde\alpha_2$, where
$\tilde\alpha_1$ is a collineation of $\PP_{n,F}$ such that
$\Gamma\tilde\alpha_1=\Gamma$, and $\tilde\alpha_2$ is a
projectivity. Since $({\rm
id}_{\displaystyle\PP_{n',F}}\times\tilde\alpha_2)\gamma=\gamma\pi$,
$\pi$ a projectivity of $\overline{\PP}$, we obtain the first
assertion with $\alpha=\tilde\alpha_1\rist{\Gamma}$ and
$\psi=\pi\tilde\psi$.

\begin{sloppypar}
Also the latter statement follows from the properties of the
action of a collineation on the Grassmann and Segre
varieties.$\Box$
\end{sloppypar}


\begin{thebibliography}{99}

\bibitem{BiHaZa99} A.~Bichara, H.~Havlicek, C.~Zanella, On linear morphisms
of product spaces, {\em Discrete Math.\/} to appear.

\bibitem{BiZa99} A.~Bichara, C.~Zanella, Characterization of embedded
special manifolds, {\em Discrete Math.\/} {\bf 208/209} (1999) 77--83.

\bibitem{Bra73} H.~Brauner, Eine geometrische Kennzeichnung linearer
Abbildungen, {\em Mh.\ Math.\/} {\bf 77} (1973) 10--20.

\bibitem{Hav81} H.~Havlicek, Zur Theorie linearer Abbildungen I, II,
{\em J.\ Geom.\/} {\bf 16} (1981) 152--180.

\bibitem{HiTh91} J. W. P.~Hirschfeld, J. A.~Thas, {\em General Galois Geometries\/},
Sect.~24.4 (Clarendon Press, Oxford, 1991).

\bibitem{Tal86} G.~Tallini, Partial line spaces and algebraic varieties,
{\em Symp.\ Math.\/} {\bf 28} (1986) 203--217.

\bibitem{Zan95} C.~Zanella, Embeddings of Grassmann spaces, {\em J. Geom.\/} {\bf
52} (1995) 193--201.

\bibitem{Zan96} C.~Zanella, Universal properties of the Corrado Segre
embedding, {\em Bull.\ Belg.\ Math.\ Soc.\ Simon Stevin\/} {\bf 3} (1996)
65--79.

\end{thebibliography}
\end{document}